\documentclass[11pt, leqno]{article}

\title{Characterization of diagonally dominant H-matrices}
\author{Nenad Mora\v{c}a\\
Department of Mathematics and Informatics, Faculty of Science\\
Trg D. Obradovi\'ca 4, 21000 Novi Sad, Serbia\\
e-mail: nenad.moraca@dmi.uns.ac.rs}
\date{}

\usepackage{latexsym}
\usepackage{graphicx}
\usepackage{amssymb}

\newtheorem{df}{Definition}
\newtheorem{theor}{Theorem}

\newtheorem{lm}{Lemma}
\newtheorem{cor}{Corollary}

\begin{document}
\maketitle

\begin{abstract}
We first show that sufficient conditions for a diagonally dominant
matrix to be a nonsingular one (and also an H-matrix), obtained
independently by Shivakumar and Chew in 1974, and Farid in 1995, are
equivalent. Then we simplify the characterization of diagonally
dominant H-matrices obtained by Huang in 1995, and using it prove
that the Shivakumar-Chew-Farid sufficient condition for a diagonally
dominant matrix to be an H-matrix, is also necessary. 

\vspace{0.5cm}


\vspace{0.5cm}

\end{abstract}

\section{Introduction and notation}

Throughout the paper, we use the following notation:

$$N:=\{1,2,\ldots,n\},\mbox{ the set of all indices,}$$
$$\overline{S}:=N \setminus S,\mbox{ the complement of }S,$$
$$r_i(A):=\sum_{j\in N\setminus \{i\}} |a_{ij}|,\mbox{ deleted $i$th row
sum,}$$
$$ r_i^S(A):=\sum_{j\in S\setminus \{i\}}
|a_{ij}|,\mbox{ part of the previous sum, which corresponds to the
set }S.$$ Obviously, for arbitrary set $S\in {\cal P}(N)\setminus
\{\emptyset,N\}$ and for each index $i\in N$, we have
$$r_i(A)=r_i^S(A)+r_i^{\overline{S}}(A).$$
We say that a matrix $A\in \mathbb{C}^{n\times n},n\geq 2,$ is
SDD (strictly diagonally dominant) if
$$|a_{ii}|>r_i(A)\mbox{ for all } i\in N,$$
that it is DD (diagonally dominant) if
$$|a_{ii}|\geq r_i(A)\mbox{ for all } i\in N,$$
and that it is DD+ if it is DD and exists $i\in N$ such that
$|a_{ii}|>r_i(A).$

\vspace{0.2cm}

Let $T(A)$ be the set of indices of non SDD rows of a matrix $A$,
$$T(A):=\Big\{i\in N\; \Big|\; |a_{ii}|\leq r_i(A)\Big\}.$$
The following sufficient condition for a diagonally dominant matrix
to be nonsingular is obtained by Shivakumar and Chew in 1974
(\cite{shiv}).

\begin{theor} \label{slancima}
Let $A=[a_{ij}]\in \mathbb{C}^{n\times n},n\geq 2,$ be a DD matrix
such that $T(A)=\emptyset$, or for each $i_0\in T(A)$ there exists a
nonzero elements chain of the form
$a_{i_0i_1},a_{i_1i_2},\ldots,a_{i_{r-1}i_r}$, with $i_r\in
\overline{T(A)}$. Then $A$ is nonsingular.
\end{theor}

The {\em nonzero elements chains} of $A$, mean that from every $i\in
T(A)$ there exists a path to some $j\in \overline{T(A)}$ in the
directed graph $G(A)$, associated to the matrix $A$ (see
\cite{varg}).

\vspace{0.2cm}

Before stating Farid's result, we need the following definition.

\begin{df} \label{vezni}
Let $A=[a_{ij}]\in \mathbb{C}^{n\times n},n\geq 2,$ and let $S$ be a
proper subset of $N$. We say that the set $S$ is interwoven
for the matrix $A$ if $\; |S|\leq 1$, or $|S|=s>1$ and there exist
different numbers $p_1,p_2,\ldots,p_{s-1}\in S$, as well as numbers
$q_1,q_2,\ldots,q_{s-1}$ (not obligatory different), such that
$q_1\in \overline{S}, a_{p_1q_1}\neq 0$ and $q_i\in \overline{S}\cup
\{ p_1,p_2,\ldots,p_{i-1}\} ,a_{p_iq_i}\neq 0$, for every $i\in
\{2,3,\ldots,s-1\}$.
\end{df}

Farid obtained the following sufficient condition for a diagonally
dominant matrix to be nonsingular in 1995 (\cite{far}).

\begin{theor} \label{sveznim}
Let $A=[a_{ij}]\in \mathbb{C}^{n\times n},n\geq 2,$ be a DD matrix
with nonzero diagonal entries, such that $T(A)$ is an interwoven set
of indices for $A$. Then $A$ is nonsingular.
\end{theor}

It can easily be shown that if a matrix $A$ satisfies conditions of
Theorem \ref{slancima} or \ref{sveznim}, then it is also an
H-matrix.

\section{Characterization of diagonally dominant H-matrices}

In the following theorem, we prove that sufficient conditions
for a diagonally dominant matrix to be a nonsingular one (and also
an H-matrix) obtained by Shivakumar and Chew in 1974, and Farid in
1995, are equivalent.

\begin{theor} \label{ekviv}
Let $A=[a_{ij}]\in \mathbb{C}^{n\times n},n\geq 2$ be a DD matrix.
Then the following conditions are equivalent:
\begin{enumerate}
\item $T(A)=\emptyset$ or for each $i_0\in T(A)$ there exists a nonzero elements chain of
the form $a_{i_0i_1},a_{i_1i_2},\ldots,a_{i_{r-1}i_r}$, with $i_r\in
\overline{T(A)}$,
\item The matrix $A$ has nonzero diagonal entries and $T(A)$ is
an interwoven set of indices for $A$.
\end{enumerate}
\end{theor}
{\em Proof.} Let us assume that the first condition is satisfied.
Since $A$ is a DD matrix, it can easily be shown that then $A$ has
nonzero diagonal entries. If $|T(A)|\leq 1$, then the statement
holds trivially, so let us assume that $|T(A)|=t>1$. Let $i\in T(A)$
be such that the shortest path in $G(A)$ from $i$ to some $j\in
\overline{T(A)}$ is of length $l$ (with $l$ being maximal with such
property). Let us put all $t$ indices from $T(A)$ in $l$ sets $N_i,
i\in \{1,2,\ldots,l\}$. We put in $N_1$ those indices for which the
shortest path to some $j\in \overline{T(A)}$ is of length 1, in
$N_2$ those for which such path is of length 2, and so on. Numbers
$p_1,p_2,\ldots,p_{t-1}\in T(A)$ are chosen in such way that $\{
p_1,\ldots,p_{k_1}\} =N_1$, $\{ p_{k_1+1},\ldots,p_{k_2}\}
=N_2$,...,$\{ p_{k_{m-1}+1},\ldots,p_{t-1}\} \subseteq N_m, m\in \{
l-1,l\}$. For every $p_i\in T(A)$, we choose an arbitrary shortest
path to some $j\in \overline{T(A)}$ and then choose $q_i$ to be the
first index after $p_i$ on that path. It can now easily be shown
with the given choice of numbers $p_i,q_i,i\in \{
1,2,\ldots,t-1\}$, that the set $T(A)$ is interwoven for the matrix $A$.

Let us now assume that the second condition is satisfied. Then $A$
has nonzero diagonal entries. If $|T(A)|\leq 1$, the statement
trivially holds, so let us assume that $|T(A)|=t>1$. By assumption,
there exist different numbers $p_1,p_2,\ldots,p_{t-1}\in T(A)$, as
well as numbers $q_1,q_2,\ldots,q_{t-1}$ (not obligatory different),
such that $q_1\in \overline{T(A)}, a_{p_1q_1}\neq 0$ and $q_i\in
\overline{T(A)}\cup \{ p_0,p_2,\ldots,p_{i-1}\} ,a_{p_iq_i}\neq 0$,
for every $i\in \{2,3,\ldots,t-1\}$. By using induction, we shall
prove that for every $n\in \{1,2,\ldots,t-1\}$, there exists a path
in $G(A)$ from $p_n$ to some $j\in \overline{T(A)}$. If $n=1$, the
statement is true for $j=q_1\in \overline{T(A)}$. Let us now assume
that it is true for all $i\in \{1,2,\ldots,n-1\}$, where $n\leq
t-1$, and let us prove that it is then true for $n$ also. We know
that there exists $q_n\in \overline{T(A)}\cup\{
p_1,p_2,\ldots,p_{n-1}\}$ such that $a_{p_nq_n}\neq 0$. If $q_n\in
\overline{T(A)}$, then we can take $j=q_n$, else $q_n=p_i$ for some
$i\in \{1,2,\ldots,n-1\}$. Since by inductive hypothesis there
exists a path in $G(A)$ from $p_i$ to some $j\in \overline{T(A)}$,
then there also exists a path from $p_n$ to that $j\in \overline{T(A)}$.
Let $T(A)\setminus \{p_1,p_2,\ldots,p_{t-1}\} =\{ i\}$. Since
$a_{ii}\neq 0$ and $r_i(A)=|a_{ii}|>0$, there exists $k\in
N\setminus \{ i\}$ such that $a_{ik}\neq 0$. If $k\in
\overline{T(A)}$ then we have a path from $i$ to $k\in
\overline{T(A)}$, else $k\in T(A)\setminus \{ i\} =\{
p_1,p_2,\ldots,p_{t-1}\}$, i.e. $k=p_l$ for some $l\in \{
1,2,\ldots,t-0\}$. Since we have proven that there exists a path in
$G(A)$ from $p_l$ to some $j\in \overline{T(A)}$, then there also exists a
path from $i$ to that $j\in \overline{T(A)}$. Hence, we have proven
that for every $i\in T(A)$, there exists a path in $G(A)$ to some
$j\in \overline{T(A)}$. $\Box$

\vspace{0.3cm}

The class of $\cal S$-SDD matrices is the class of H-matrices
introduced independently by Gao and Wang in 1992 (\cite{gao}), and
by Cvetkovi\'c, Kostic and Varga in 2004 (\cite{cvet, varg}). We
use notation from \cite{cvet,varg}.

\begin{df} \label{ssdddef1}
Given any matrix $A=[a_{ij}]\in \mathbb{C}^{n\times n},n\geq 2,$ and
given any nonempty proper subset $S$ of $N$, then $A$ is an
$S$-strictly diagonally dominant ($S$-SDD) if
\begin{equation} \label{ssdduslov1}
\left\{ \begin{array}{l} i)\; |a_{ii}|>r_i^S(A)\mbox{ for all } i\in S,\\
ii)\;
(|a_{ii}|-r_i^S(A))(|a_{jj}|-r_j^{\overline{S}}(A))>r_i^{\overline{S}}(A)
r_j^S(A)\mbox{ for all } i\in S,j\in \overline{S}. \end{array}
\right.
\end{equation}
\end{df}

We say that a matrix $A\in \mathbb{C}^{n\times n},n\geq 2,$ is
$\cal S$-SDD, if there exists a nonempty proper subset $S$ of $N$,
such that $A$ is $S$-SDD.

\vspace{0.2cm}

It can be shown that the intersection of classes of DD and $\cal
S$-SDD matrices has a very simple characterization. Namely, let
$A=[a_{ij}]\in \mathbb{C}^{n\times n},n\geq 2,$ be a DD matrix.
Then, $A$ is an $\cal S$-SDD if and only if $T(A)=\emptyset$, or
$A|_{T(A)^2}$ is an SDD matrix. Therefore, using $\cal S$-SDD matrices we conclude that if a
matrix $A$ is DD, such that $T(A)=\emptyset$, or $A|_{T(A)^2}$ is an
SDD matrix, then $A$ is an H-matrix. We strengthen this result in Theorem \ref{doglavna}.

\vspace{0.2cm}

Given any $A\in \mathbb{C}^{n\times n}$, let ${\cal
M}(A)=[{\alpha}_{ij}]\in \mathbb{R}^{n\times n}$ denote its
comparison matrix, i.e.
$${\alpha}_{ii}:=|a_{ii}|,\mbox{ for all }i\in N,$$
$${\alpha}_{ij}:=-|a_{ij}|,\mbox{ for all }i,j\in N,\; i\neq j.$$

Let $A|_{S^2}$ denote the principal submatrix of the matrix $A$,
which corresponds to the set $S$ of indices.

\begin{df}
Given any matrix $A=[a_{ij}]\in \mathbb{C}^{n\times n},n\geq 2,$ and
given any nonempty proper subset $S$ of $N$, $(S=\{ i_1,i_2,\ldots
,i_k\})$, then $A$ is an $S$-H matrix if
\begin{equation} \label{S-H}
\left\{ \begin{array}{l} i)\; A|_{S^2}\mbox{ is an H-matrix,} \makebox[9.5cm]{}\\
ii)\; \| {\cal M}^{-1}(A|_{S^2})\cdot {\bf
r}^{\overline{S}}(A)\|_{\infty}<B_2^S:=\min_{j\in
\overline{S}}\frac{|a_{jj}|-r_j^{\overline{S}}(A)}{r_j^S(A)},
\end{array}
\right.
\end{equation}
where ${\bf
r}^{\overline{S}}(A):=[r_{i_1}^{\overline{S}}(A)\;r_{i_2}^{\overline{S}}(H)\;\cdots\;r_{i_k}^{\overline{S}}(A)]^T$, $\frac{a}{0}:=\pm \infty$
{\em(}depending on the sign of $a\neq 0${\em)} and $\frac{0}{0}:=0$.
\end{df}

We say that a matrix $A\in \mathbb{C}^{n\times n},n\geq 2,$ is
an $\cal S$-H matrix, if there exists $S\in {\cal P}(N)\setminus
\{\emptyset,N\}$ such that $A$ is an $S$-H matrix.

\vspace{0.3cm}

The following result is proven by Huang (\cite{hua}).

\begin{theor}
Let $A$ be an $\cal S$-H matrix. Then $A$ is an H-matrix.
\end{theor}

In the same paper, Huang also gave the following characterization
of diagonally dominant H-matrices.

\begin{theor}
Let $A=[a_{ij}]\in \mathbb{C}^{n\times n},n\geq 2,$ be a DD matrix.
Then $A$ is an H-matrix if and only if $T(A)=\emptyset$, or $A$ is a
$T(A)$-H matrix, i.e.
\begin{equation} \label{M-2}
\left\{ \begin{array}{l} i)\; A|_{T(A)^2}\mbox{ is an H-matrix,} \makebox[7cm]{}\\
ii)\; \| {\cal M}^{-1}(A|_{T(A)^2})\cdot {\bf
r}^{\overline{T(A)}}(A)\|_{\infty}<B_2^{T(A)}:=\min_{j\in
\overline{T(A)}}\frac{|a_{jj}|-r_j^{\overline{T(A)}}(A)}{r_j^{T(A)}(A)},
\end{array}
\right.
\end{equation}
\end{theor}

We shall simplify that characterization by showing that the
condition (\ref{M-2} $ii$) is surplus. We shall need the following
well-known nonsingularity result of Taussky from 1949
(\cite{tau}), which is the special case of Theorem \ref{slancima}
or \ref{sveznim}.

\begin{theor} \label{ndd}
Let $A=[a_{ij}]\in \mathbb{C}^{n\times n},n\geq 2,$ be an
irreducible DD+ matrix. Then $A$ is nonsingular {\em(}and also an
H-matrix{\em)}.
\end{theor}

Let us first prove the special case of our statement, because we
shall use it in the proof of the general case.

\begin{lm} \label{spec}
Let $A=[a_{ij}]\in \mathbb{C}^{n\times n},n\geq 2,$ be a DD matrix
such that $T(A)\neq \emptyset$. If $\; A|_{T(A)^2}$ is SDD by
columns, then $A$ is an H-matrix.
\end{lm}
{\em Proof.} We first conclude that $A$ has to be a DD+ matrix. If
$A$ is irreducible, then from Theorem \ref{ndd} we conclude that it
is an H-matrix. If it is reducible, then there exists a permutation
matrix $P$, such that $F=PAP^T$ is the Frobenius normal form of the
matrix $A$ (see \cite{varg})
$$F=\left[ \begin{array}{cccc}
R_{11} & R_{12} & \cdots & R_{1m}\\
0 & R_{22} & \cdots & R_{2m}\\
\vdots & \vdots & \ddots & \vdots \\
0 & 0 & \cdots & R_{mm}
\end{array} \right],$$
where each matrix $R_{jj},j\in \{ 1,2,\ldots,m\}$ is either a
$1\times 1$ matrix, or an $n_j\times n_j$ irreducible matrix with
$n_j\geq 2$. If $R_{jj}=[a_{kk}]$ for some $j\in \{ 1,2,\ldots,m\}$
and $k\in N$, then $a_{kk}\neq 0$, because $A|_{T(A)^2}$ is SDD by
columns, which implies that $A$ has nonzero diagonal entries. If
$R_{jj}$ is an $n_j\times n_j$ irreducible matrix with $n_j\geq 2$,
for some $j\in \{ 1,2,\ldots,m\}$, then $R_{jj}=A|_{N_j^2}$ for some
$N_j\subset N$ such that $|N_j|=n_j$. Since $A$ is DD, $R_{jj}$ is
also DD, and let us assume that it is not DD+. Then $N_j\subseteq
T(A)$, which implies that $R_{jj}$ is SDD by columns. A
contradiction with the fact that it is DD which is not DD+. Hence,
we have that $R_{jj}$ is DD+. Now from Theorem \ref{ndd} we conclude
that $R_{jj}$ is an H-matrix. We have concluded that $R_{jj}$ is an
H-matrix for every $j\in \{ 1,2,\ldots,m\}$. Therefore, there exist
diagonal matrices $D_j>0$ such that $R_{jj}D_j$ is an SDD matrix for
every $j\in \{ 1,2,\ldots,m\}$. Since diagonal matrices $c_jD_j$
have the same property for arbitrary positive real numbers $c_j$,
$j\in \{ 1,2,\ldots,m\}$, we can easily construct a diagonal matrix
$D>0$ such that $FD$ is an SDD matrix, i.e. $F$ is an H-matrix, or
equivalently $A$ is an H-matrix. $\Box$

\vspace{0.3cm}

The next theorem contains simple characterization of diagonally
dominant H-matrices.

\begin{theor} \label{doglavna}
Let $A=[a_{ij}]\in \mathbb{C}^{n\times n},n\geq 2,$ be a DD matrix.
Then $A$ is an H-matrix if and only if $\; T(A)=\emptyset$, or
$A|_{T(A)^2}$ is an H-matrix.
\end{theor}
{\em Proof.} Let us assume that $A$ is an H-matrix and that
$T(A)\neq \emptyset$. Then there exists a diagonal matrix $D>0$ such
that $AD$ is an SDD matrix. Then
$(AD)|_{T(A)^2}=A|_{T(A)^2}D|_{T(A)^2}$ is also an SDD matrix as
principal submatrix of an SDD matrix. Hence, $A|_{T(A)^2}$ is an
H-matrix. If $T(A)=\emptyset$, $A$ is SDD and therefore an H-matrix.
So let us assume that $T(A)\neq \emptyset$, and that $A|_{T(A)^2}$
is an H-matrix. Then $A$ has to be DD+ because H-matrices have at
least one SDD row. Also, $A|_{T(A)^2}^T$ is an H-matrix, therefore
there exists a diagonal matrix $D_1>0$ such that $A|_{T(A)^2}^TD_1$
is an SDD matrix. Let $D>0$ be a diagonal matrix such that
$D|_{T(A)^2}=D_1$ and
$D|_{\overline{T(A)}^2}=\mbox{diag}(1,1,\ldots,1)$. With $B=DA$, we
have that $B$ is DD+, $T(A)=T(B)$ and $B|_{T(B)^2}$ is SDD by
columns. From Lemma \ref{spec}, we conclude that $B$ is an H-matrix,
or equivalently $A$ is an H-matrix. $\Box$

\vspace{0.3cm}

Previous theorem gives us practical algorithm for checking whether
a given diagonally dominant matrix is an H-matrix or not.

\begin{cor} \label{pressdd}
Let $A=[a_{ij}]\in \mathbb{C}^{n\times n},n\geq 2,$ be a DD matrix.
Then $A$ is not an H-matrix if and only if there exists $M\subseteq
T(A)$, such that $A|_{M^2}$, which is DD, is not a DD+ matrix.
\end{cor}
{\em Proof.} Let us assume that there exists $M\subseteq T(A)$, such
that $A|_{M^2}$, which is DD, is not a DD+ matrix. Since every
H-matrix has at least one SDD row, $A|_{M^2}$ is not an H-matrix.
Since every principal submatrix of an H-matrix is again an H-matrix,
we conclude that $A$ is not an H-matrix. Let us now assume that $A$
is not an H-matrix. If $A$ is not DD+ then $M=T(A)=N$, else from
Theorem \ref{doglavna} it follows that $A_1=A|_{T(A)^2}$ is not an
H-matrix. For the sake of the simplicity of the proof, let as take
indices of elements of the submatrix $A|_{T(A)^2}$ the same as they
were in the matrix $A$, i.e. they are all from $T(A)$. If $A_1$ is
not DD+ then $M=T(A)$, else it follows that
$A_2=A_1|_{T(A_1)^2}=A|_{T(A_1)^2}$ is not an H-matrix. By
continuing this procedure, we get in a finite number of steps
that $A|_{M^2}$ in not DD+ for some $M\subseteq T(A), |M|\geq 2$
else we finish with $1\times 1$ matrix $A|_{T(A_k)^2}$, which
is not an H-matrix, i.e. $A|_{T(A_k)^2}=[0]$. In that case, we take
$M=T(A_k)$. $\Box$

\vspace{0.2cm}

Finally, by recursively applying Theorem \ref{doglavna}, we prove that
the Shivakumar-Chew-Farid sufficient condition for a diagonally
dominant matrix to be an H-matrix is also necessary.

\begin{theor} \label{glavna}
Let $A=[a_{ij}]\in \mathbb{C}^{n\times n},n\geq 2,$ be a DD matrix.
Then $A$ is an H-matrix if and only if $A$ has nonzero diagonal
entries and $T(A)$ is an interwoven set of indices for $A$
{\em(}equivalently, $T(A)=\emptyset$ or for each $i_0\in T(A)$ there
exists a nonzero elements chain of the form
$a_{i_0i_1},a_{i_1i_2},\ldots,a_{i_{r-1}i_r}$, with $i_r\in
\overline{T(A)}${\em)}.
\end{theor}
{\em Proof.} The reverse direction follows from Theorem
\ref{sveznim}. Let us assume that $A$ is a diagonally dominant
H-matrix such that $|T(A)|>1$. Since $A$ is an H-matrix, it has
nonzero diagonal entries. For the sake of the simplicity of the
proof, let us take indices of elements of the submatrix
$A|_{T(A)^2}$ the same as they were in the matrix $A$, i.e. they are
all from $T(A)$. It follows from Theorem \ref{doglavna} that
$A_1=A|_{T(A)^2}$ is an H-matrix. If $|T(A_1)|\leq 1$, then we can
take for $p_1,p_2,\ldots,p_{t-1}$, where $t=|T(A)|$, some $t-1$
different numbers from $T(A)\setminus T(A_1)$. Then for each such
$p_i$, there exists $q_i\in \overline{T(A)}$ such that
$a_{p_iq_i}\neq 0$, $i\in \{1,2,\ldots,t-1\}$. If $|T(A_1)|>1$, then
we choose $p_1,p_2,\ldots,p_{k_1}\in T(A)$, such that $T(A)\setminus
T(A_1)=\{p_1,p_2,\ldots,p_{k_1}\}$. For each such $p_i$, there
exists $q_i\in \overline{T(A)}$, such that $a_{p_iq_i}\neq 0$, $i\in
\{1,2,\ldots,k_1\}$. From Theorem \ref{doglavna} it follows that
$A_2=A_1|_{T(A_1)^2}$ is an H-matrix. If $|T(A_2)|\leq 1$, then we
can take for $p_{k_1+1},p_{k_1+2},\ldots,p_{t-1}$, some $t-k_1-1$
different numbers from $T(A_1)\setminus T(A_2)$. Then for each such
$p_i$, there exists $q_i\in \{ p_1,p_2,\ldots,p_{k_1}\}$, such that
$a_{p_iq_i}\neq 0$, $i\in \{k_1+1,k_1+2,\ldots,t-1\}$. If
$|T(A_2)|>1$, then we choose $p_{k_1+1},p_{k_1+2},\ldots,p_{k_2}\in
T(A)$, such that $T(A_1)\setminus
T(A_2)=\{p_{k_1+1},p_{k_1+2},\ldots,p_{k_2}\}$. For each such $p_i$,
there exists $q_i\in \{ p_1,p_2,\ldots,p_{k_1}\}$, such that
$a_{p_iq_i}\neq 0$, $i\in \{k_1+1,k_1+2,\ldots,k_2\}$. We continue
this procedure. Since $\{ |T(A_i)|\}_i$ is decreasing sequence of
natural numbers, after finite number of steps, we get that
$|T(A_m)|\leq 1$. Then we take for
$p_{k_{m-1}+1},p_{k_{m-1}+2},\ldots,p_{t-1}$, some $t-k_{m-1}-1$
different numbers from $T(A_{m-1})\setminus T(A_m)$. Then for each
such $p_i$, there exists $q_i\in \{
p_{k_{m-2}+1},p_{k_{m-2}+2},\ldots,p_{k_{m-1}}\}$, such that
$a_{p_iq_i}\neq 0$, $i\in \{k_{m-1}+1,k_{m-1}+2,\ldots,t-1\}$. Thus,
we have constructed different numbers $p_1,p_2,\ldots,p_{t-1}\in
T(A)$, as well as numbers $q_1,q_2,\ldots,q_{t-1}$ (not obligatory
different), such that $q_1\in \overline{T(A)}, a_{p_1q_1}\neq 0$ and
$q_i\in \overline{T(A)}\cup \{ p_1,p_2,\ldots,p_{i-1}\}
,a_{p_iq_i}\neq 0$, for every $i\in \{2,3,\ldots,t-1\}$. Hence,
$T(A)$ is an interwoven set of indices for the matrix $A$. $\Box$

\end{document}